\documentclass[12pt]{amsart}
\usepackage{amsmath}
\usepackage{graphicx,color,epsfig}
\theoremstyle{plain}
\newtheorem{theorem}                 {Theorem}      [section]
\newtheorem{proposition}  [theorem]  {Proposition}

\theoremstyle{definition}

\newtheorem{remark}       [theorem]  {Remark}
\newtheorem{definition}   [theorem]  {Definition}

\numberwithin{equation}{section} \numberwithin{theorem}{section}
\title[SOME CLASSES OF GEODESICS ]{GEODESICITY AND ISOCLINITY PROPERTIES FOR THE TANGENT BUNDLE
OF THE HEISENBERG MANIFOLD WITH SASAKI METRIC}
\author{S. L. DRU\c T\u A and M. P. PIU}
\date{}
\begin{document}
\maketitle
\begin{abstract}
We prove that the horizontal and vertical distributions of the
tangent bundle with the Sasaki metric are isocline, the
distributions given by the kernels of the horizontal and vertical
lifts of the contact form $\omega$ from the Heisenberg manifold
$(H_3,g)$ to $(TH_3,g^S)$ are not totally geodesic, and the
distributions $F^H=L(E_1^H,E_2^H)$ and $F^V=L(E_1^V,E_2^V)$ are
totally geodesic, but they are not isocline. We obtain that the
horizontal and natural lifts of the curves from the Heisenberg
manifold $(H_3,g)$, are geodesics in the tangent bundle endowed
with the Sasaki metric $(TH_3,g^s)$, if and only if the curves
considered on the base manifold are geodesics. Then, we get two
particular examples of geodesics from $(TH_3,g^s)$, which are not
horizontal or natural lifts of geodesics from the base manifold
$(H_3,g)$.
\end{abstract}

\vskip2mm
{\bf Mathematics Subject Classification 2000:}


\section{Introduction}

The tangent bundle $TM$ of a Riemannian manifold splits into the
vertical and horizontal distributions, defined by the Levi Civita
connection of the metric $g$ from the base manifold (see
\cite{YanoIsh}).

In the study of the differential geometry of the tangent bundle of
a Rieman-\ nian manifold, one uses several (pseudo) Riemannian
metrics, induced by the Riemannian metric from the base manifold,
and constructed on the horizontal and vertical distributions.

Maybe the best known Riemannian metric on the tangent bundle is
that introduced by Sasaki in the paper \cite{Sasaki} from 1958.
The results form \cite{Janyska}--\cite{KowalskiSek}, concerning
the natural lifts, allowed the extension of the Sasaki metric, to
the metrics of natural diagonal lift type (see \cite{Oproiu3}) and
general natural lifted metrics (see \cite{Oproiu4},
\cite{Tahara2}), leading to interesting geometric structures
studied in the last years (see \cite{Abbassi1}, \cite{Munteanu1}
-- \cite{OprPap1}), and to interesting relations with some
problems in Lagrangian and Hamiltonian mechanics (see
\cite{Anast}, \cite{Miron}, \cite{MironAnast}).

The aim of the first section from this paper is to find some
geometric pro-\ perties of the horizontal and vertical
distributions of the tangent bundle $TM$ of a Riemannian manifold
$(M,g)$, endowed with the Sasaki metric $g^s$. More precisely we
shall study the property of the the two distributions of being
isocline (and implicitly totally geodesic), with respect to the
Sasaki metric.

The notion of \emph{isocline distribution} was introduced by Lutz,
which made in \cite{lutz} a metric study of the contact
structures, measuring, with the help of a metric $g$, the
evolution of a contact structure $F$ along the geodesics of $g$.
One of the metric characters of a totally geodesic field $F$,
considered by Lutz in \cite{lutz}, is the evolution of its angle
along an arbitrary geodesic $\gamma$. When the angle between a
totally geodesic distribution $F$ and $\dot\gamma$ is constant
along the geodesic $\gamma$, the totally geodesic distribution is
called \emph{isocline}.

The second author studied in her PhD thesis \cite{Piu}, the
property of being isocline for the contact structures on
hyper-surfaces of $\mathbb{R}^{2n+2}$. In the third section of the
present paper we shall prove that the horizontal and the vertical
distribution of the tangent bundle of a Riemannian manifold are
always isocline with respect to the Sasaki metric $g^s$, and we
shall construct some examples of distributions on the tangent
bundle of the Heisenberg manifold, which have or have not the
properties of being totally geodesic or isocline, respectively. To
this aim, we shall consider the horizontal and vertical lifts of
the contact form $\omega$ from the Heisenberg manifold $(H_3,g)$
to the tangent bundle $(TH_3,g^s)$, and we shall prove that the
distributions given by $F=Ker(\omega^H)$ (or by $F=Ker(\omega^V)$)
are not totally geodesic, but the distributions
$F^H=L(E_1^H,E_2^H)$ and $F^V=L(E_1^V,E_2^V)$ are totally geodesic
and not isocline.

An important geometric problem is to find the geodesics on the
smooth manifolds with respect to the Riemannian metrics (see
\cite{Boothby}--\cite{Gheorghiev}, \cite{Nagy},
\cite{Piu2}-\cite{SalimovK}, \cite{YanoIsh}). In \cite{YanoIsh},
Yano and Ishihara proved that the curves on the tangent bundles of
Riemannian manifolds are geodesics with respect to certain lifts
of the metric from the base manifold, if and only if the curves
are obtained as certain types of lifts of the geodesics from the
base manifold. In two very recent papers, Salimov and his
collaborators studied the analogous problem for the geodesics on
the tangent bundles endowed with Cheeger-Gromoll metrics (see
\cite{SalimovK}), and on the tensor bundles with Sasakian metrics
(see \cite{Salimov}) .

In the last section of the present paper, we are interested in
finding some concrete examples of geodesics on the tangent bundle
$(TH_3,g^s)$, of the Heisenberg manifold $(H_3,g)$,  with respect
to the Sasaki metric $g^s$.

We prove that if $C$ is a curve in the Heisenberg manifold
$(H_3,g)$, then its horizontal and natural lifts, $\widetilde{C}$
and $\widehat{C}$, passing through the origin, such that
$\dot{\widetilde{C}}(0)=\dot{\widehat{C}}(0)=(u,v,w,0,0,0)$, are
geodesics on $(TH_3,g^s)$, if and only if the curve considered on
the base manifold is a geodesic.

Working in a more general context, we look for some classes of
geodesics on $(TH_3,g^s)$, which are not obtained as horizontal or
natural lifts of the geodesics from the base manifold.

\section{Preliminary results.} Let $(M,g)$ be a smooth $n$-dimensional Riemannian manifold and
denote its tangent bundle by $\tau :TM\longrightarrow M$. Just to
fix the notation, recall some basic things about $TM$. It has a
structure of a $2n$-dimensional smooth manifold, induced from the
smooth manifold structure of $M$. This structure is obtained by
using local charts on $TM$ induced  from usual local charts on
$M$. If $(U,\varphi )= (U,x^1,\dots ,x^n)$ is a local chart on
$M$, then the corresponding induced local chart on $TM$ is $(\tau
^{-1}(U),\Phi )=(\tau ^{-1}(U),x^1,\dots , x^n,$ $y^1,\dots ,y^n)$
(see \cite{YanoIsh} for further details).

Denote by $\nabla$ the Levi Civita connection of the Riemannian
metric $g$ on $M$. Then we have the direct sum decomposition
\begin{equation*}\label{decomp}
TTM=VTM\oplus HTM
\end{equation*}
of the tangent bundle to $TM$ into the vertical distribution
$VTM={\rm Ker}\ \tau_*$ and the horizontal distribution $HTM$
defined by $\nabla $.

The set set of vector fields $\big\{\frac{\partial}{\partial
y^i},\frac{\delta}{\delta x^j}\big\}_{i,j=\overline{1,n}}$ defines
a local frame on $TM$, adapted to the direct sum decomposition
\eqref{decomp}. Notice that
$$
\frac{\partial}{\partial y^i}=\left(\frac{\partial}{\partial
x^i}\right)^V,\ \ \frac{\delta}{\delta
x^j}=\left(\frac{\partial}{\partial
x^j}\right)^H=\frac{\partial}{\partial x^j}-\Gamma^h_{0j}
\frac{\partial}{\partial y^h},\ \ \ \Gamma ^h_{0j}=y^k\Gamma
^h_{kj}
$$
where $X^V\in VTM$ and $X^H\in HTM$ denote the vertical and
horizontal lift of the vector field $X$ on $M$, respectively, and
$\Gamma ^h_{kj}(x)$ are the Christoffel symbols of $g$.

The Sasaki metric $g^s$ on the tangent bundle $TM$ is defined by
the relations
$$
\begin{cases}
g^S(X^H,Y^H)=g^S(X^V,Y^V)=g(X,Y)\circ \tau\\
g^S(X^H,Y^V)=g^S(X^V,Y^H)=0
\end{cases}
\forall X,Y\in \mathcal{T}^1_0(M).
$$

If the metric $g$ from the base manifold $M$ has the components
$g_{ij}$ in a coordinate neighborhood, then the Sasaki metric
$g^s$ on the tangent bundle may be defined as the Riemannian
metric which has the expression
\begin{equation*}
g^s=g_{ij}dx^idx^j+g_{ij}Dy^iDy^j,\quad \forall i,j=\overline{1,n}
\end{equation*}
where $\big\{Dy^i,dx^j\big\}_{i,j=\overline{1,n}}$ is the dual
frame of $\big\{\frac{\partial}{\partial y^i},\frac{\delta}{\delta
x^j}\big\}_{i,j=\overline{1,n}}$. The covariant derivative of
$y^i$ with respect to the Levi-Civita connection of the metric $g$
is given by
\begin{equation*}\label{Dyi}
Dy^i=dy^i+\Gamma^i_{0j}dx^j,\qquad \Gamma^i_{0j}=\Gamma^i_{hj}y^h.
\end{equation*}

In particular, if the base manifold is the Heisenberg manifold
$(H_3,g)$, where
\begin{equation}\label{Heisenberg}
g=(dx^1)^2+(dx^2)^2+(dx^3+x^2dx^1-x^1dx^2)^2
\end{equation}
then the Sasaki metric on $TH_3$ is
\begin{align}\label{Sasaki}
g^s&=(dx^1)^2+(dx^2)^2+(dx^3+x^2dx^1-x^1dx^2)^2
\\
&+(Dy^1)^2+(Dy^2)^2+(Dy^3+x^2Dy^1-x^1Dy^2)^2.\nonumber
\end{align}

\section{Totally geodesic and isocline distributions on the tangent bundle of a Riemannian manifold}

It is well known that the horizontal distribution $HTM$ of the
tangent bundle of  an $n$-dimensional Riemannian
manifold $(M,g)$, is integrable if and only if the manifold $M$ is
flat. In this section we shall prove that the horizontal and the
vertical distributions of $TM$ are always isocline with respect to
the Sasaki metric $g^s$.

A regular distribution $F$ defined on a connected Riemannian
manifold $(M,g)$ is called totally geodesic if every geodesic
tangent to $F$ in one point is tangent to the distribution in all
the points.

A distribution $F$ is totally geodesic (\cite{Wood}) if and only
if the distribution $D$ normal to $F$ is Riemannian (i.e. if ($L_Z
g)(X,Y) = 0$ for every $X,Y \in C^\infty(F)$ and $Z \in
C^\infty(D)$).

A totally geodesic distribution forms a constant angle  with the
integral geodesic curve. If, more over, the field of planes makes
a constant angle with the tangent vector field $\dot{\gamma}$
along an arbitrary geodesic curve $\gamma$ we say that the
structure is isocline.

Let $F$ be a totally geodesic distribution and $N$ a unitary vector
field, normal to $F$.

\begin{definition}
The totally geodesic contact structure $F$ is called isocline if
for every geodesic curve parameterized by the arc length, the
angle between $F$ and the tangent vector field $\dot{\gamma}(s)$
is constant along the geodesic.
\end{definition}

\begin{proposition}\cite{lutz}
If $\nabla$ is the Levi-Civita connection associated to the metric
$g$ on $M$, a totally geodesic  distribution $F$ is isocline if and
only if for every vector field $N$ normal to $F$, the vector field
$ {\nabla}_NN$ is normal to $F$.
\end{proposition}
If $\{X_i,N_\alpha\}$, $i,=\overline{1,p},\alpha=\overline{1,q},p
+ q = n,$ is an orthonormal frame of $(M,g)$ adapted to the
distribution $F$ ($X_i \in C^\infty(F)$ and $N_\alpha \in
C^\infty(F^\bot)$) then $F$ is isocline if and only if
\begin{eqnarray}\label{t.g}
g\left({\nabla}_{X_i}X_j +{\nabla}_{X_j}X_i, N_\alpha\right) & = &0 \qquad  \mbox{geodesicity}\\
g\left(\nabla_{N_\alpha}N_\beta + \nabla_{N_\beta}N_\alpha
,X_i\right)  &=& 0 \label{isocline}
\end{eqnarray}
where $i, j=\overline{1,p}; \alpha, \beta=\overline{1,q}$.

\begin{proposition}
The distributions $HTM$ and $VTM$ are isocline.
\end{proposition}
\begin{proof}
Let us consider the connection $\bar{\nabla}$ of $(TM,g^s)$ and
the Levi-Civita connection $\nabla$ of $(M,g)$, which are related
by the formulas (see \cite{Blair})
\begin{eqnarray*}\label{LC}
\begin{cases}
(\bar{\nabla}_{X^H} Y^H)_Z = (\nabla_X Y)^H_Z -
\frac{1}{2} (R_{XY}Z)^V \\
(\bar{\nabla}_{X^H} Y^V)_Z = (\nabla_X Y)^V_Z -
\frac{1}{2} (R_{YZ}X)^H\\
(\bar{\nabla}_{X^V} Y^H)_Z =  -
\frac{1}{2} (R_{XY}Z)^V\\
(\bar{\nabla}_{X^V} Y^V) = 0
\end{cases}
\end{eqnarray*}
where $X,Y$ are tangent to $M$.

Now, we may easily verify the conditions (\ref{t.g})  and
(\ref{isocline}) for $HTM$ and $VTM$ to be totally geodesic and
isocline:
\begin{eqnarray*}
g^s(\bar{\nabla}_{X^H} Y^H + \bar{\nabla}_{Y^H} X^H, X^V) &=& 0 \\
g^s(\bar{\nabla}_{X^V} Y^V + \bar{\nabla}_{Y^V} X^V, Y^H) &=& 0.
\end{eqnarray*}
\end{proof}

In the sequel, we shall focus our attention to the geometry of the
Heisenberg manifold, usually known as Heisenberg group $H_3$. A
first remark is that its contact distribution furnishes an example
of totally geodesic and isocline distribution (see \cite{Hangan},
\cite{lutz}).

We are interested in finding examples of distributions on the
tangent bundle of the Heisenberg manifold, $TH_3$, endowed with
the Sasaki metric expressed by (\ref{Sasaki}), which are isocline,
or which are only totally geodesic, without being isocline.

The metric $g$ on $H_3$, given by (\ref{Heisenberg}), is invariant
with respect to the left translations and with respect to the
rotations around the $z$ axis. We shall use the invariant
orthonormal coframe
$$\theta^1=dx^1, \quad \theta^2=dx^2, \quad \theta^3 = dx^3+x^2dx^1-x^1dx^2$$
and the dual basis
$$E_1=\frac{\partial}{\partial x^1}- x^2\frac{\partial}{\partial x^3}, \quad
E_2=\frac{\partial}{\partial x^2}+x^1\frac{\partial}{\partial x^3}, \quad
E_3=\frac{\partial}{\partial x^3}.$$

The Levi-Civita conection of the metric $g$ is given by
\begin{equation*}
\left\{
\begin{array}{l}
\nabla_{E_1}E_1=0, \quad \nabla_{E_1}E_2 = E_3, \quad
\nabla_{E_1}E_3=- E_2
\\ \mbox{} \\
\nabla_{E_2}E_1=- E_3, \quad \nabla_{E_2}E_2=0, \quad
\nabla_{E_2}E_3= E_1
\\ \mbox{} \\
\nabla_{E_3}E_1=- E_2, \quad \nabla_{E_3}E_2 = E_1, \quad
\nabla_{E_3}E_3=0.
\end{array}
\right.
\label{eq:conn_1}
\end{equation*}

The non vanishing components of the curvature tensor field
$$
R(X,Y)Z = \nabla_X\nabla_YZ - \nabla_Y\nabla_XZ-\nabla_{[X,Y]}Z
$$
and of the Riemann-Christoffel curvature $R(X,Y,Z,W) =
g(R(X,Y)W,Z)$  are
\begin{equation*}
\left\{
\begin{array}{l}
R^2_{112}= 3, \quad R^3_{113}=- 1, \quad R^1_{212}= - 3
 \\ \mbox{} \\
R^1_{313}= 1,\quad R^3_{223}= - 1, \quad R^2_{323}=  1
\end{array}
\right.
\label{eq:curb_1}
\end{equation*}
\begin{equation*}
R_{1212}= -3, \quad R_{1313} = R_{2323}= 1
\end{equation*}
where we used the notations
$$
R(E_a,E_b)E_c = R^i_{cab}E_i, \quad R(E_a,E_b,E_c,E_d) = R_{abcd}.
$$

We may ask what happens with the distributions determined by the
kernels of the horizontal and vertical lift  of the contact form
$\omega = dx^3 + x^2 dx^1 - x^1 dx^2$. In this sense, we may prove
the following proposition.

\begin{proposition}
If $\omega^H$ ($\omega^V$) is the horizontal (vertical) lift of
the contact form $\omega$ from the Heisenberg manifold $H_3$, then
the distribution $F$ of codimension $1$, defined by $F =
Ker(\omega^H)$ ($F = Ker(\omega^V)$) is not totally geodesic.
\end{proposition}
\begin{proof}
If the distribution $F$ is given by $Ker(\omega^H)$, then one can
choose a basis given by the vector fields
$\{E_1^H,E_2^H,E_1^V,E_2^V,E_3^V\}$, and we may easily verify that
\begin{eqnarray*}
g^s_{(p,y)}\left(\widetilde{\nabla}_{E_1^H} E_1^V +
\widetilde{\nabla}_{E_1^V }
E_1^H, E_3^H\right) &=& - \frac{1}{2}g^s_{(p,y)}\left((R(E_1,y)E_1)^H,E_3^H\right) \\
&=& -\frac{1}{2}g_{\tau(y)}(R(E_1,y)E_1,E_3) \neq 0 \nonumber
\end{eqnarray*}
where $\widetilde{\nabla}$ is the Levi-Civita connection of the
Sasaki metric on $TH_3$, and $y$ represents a tangent vector from
$TH_3$.

Analogously, if the distribution is defined by $F =
Ker(\omega^V)$, then a basis for $F$ is given by the vector fields
$\{E_1^H,E_2^H,E_3^H,E_1^V,E_2^V\}$, and it may be verified that
\begin{eqnarray*}
 g^s_{(p,y)}\big(\widetilde{\nabla}_{E_1^H}
E_2^V + \widetilde{\nabla}_{E_2^V } E_1^H, E_3^V\big) &=&
g^s_{(p,y)}\big((\nabla_{E_1} E_2)^V, E_3^V\big)-
 \\
&-&\frac{1}{2}g^s_{(p,y)}\big((R(E_2,y)E_1)^H-
(R(E_2,E_1)y)^V,E_3^V\big)\\
& =& g_{\tau(y)}(\nabla_{E_1} E_2,E_3) = 1
\neq 0.
\end{eqnarray*}

Thus the proposition is proved.
\end{proof}

Now we give an example of totally geodesic distributions on
$TH_3$, which are not isocline.

\begin{proposition}
The distributions $F^H = L(E_1^H,E_2^H)$ and $F^V =
L(E_1^V,E_2^V)$ are totally geodesic and they are
not isocline.
\end{proposition}
\begin{proof}
Since the Levi-Civita connection of the Sasaki metric $g^s$ from
the tangent bundle of a Riemannian manifold $(M,g)$ has the
expressions (\ref{LC}), and the Levi-Civita connection from
$(H_3,g)$ satisfies the relation (\ref{eq:conn_1}), we may easily
prove that the Levi-Civita connection $\widetilde{\nabla}$ from
$TH_3$ verify the relations
\begin{eqnarray*}
\widetilde{\nabla}_{E_1^H} E_2^H + \widetilde{\nabla}_{E_2^H} E_1^H  &=& 0 \\
\widetilde{\nabla}_{E_1^V} E_2^V + \widetilde{\nabla}_{E_2^V}
E_1^V &=& 0
\end{eqnarray*}
and thus both distributions $F^H,F^V$ are totally geodesic.

We may easily prove that $\widetilde{\nabla}$ fulfills also the
relations
\begin{eqnarray*}
 g^s_{(p,y)}\left(\widetilde{\nabla}_{E_3^H} E_1^V +
 \widetilde{\nabla}_{E_1^V} E_3^H , E_1^H\right)
 &=&-\frac{1}{2} g^s_{(p,y)}\left(\big(R(E_1,y)E_3\big)^H,E_1^H\right)  \neq 0 \\
 g^s_{(p,y)}\left(\widetilde{\nabla}_{E_3^V} E_1^H +
 \widetilde{\nabla}_{E_1^H} E_3^V, E_1^V\right)&=&-\frac{1}{2}
 g^s_{(p,y)}\left(\big(R(E_3,E_1)y\big)^V,E_1^V\right)  \neq 0 \end{eqnarray*}
which prove that the distributions $F^H$ and $F^V$ are not
isocline.
\end{proof}

\section{Geodesics in the tangent bundle $\mathbf{(TH_3,g^s)}$.}

Let $M$ be an $n$-dimensional Riemannian manifold and $C:I \rightarrow M$ a curve
parametrized on it, expressed locally by
$$
C(t) = \{x^1(t), \cdots, x^n(t)\},
$$
and let $X$ be a vector field along a curve $C$. Then, in the
tangent bundle $TM$, a curve $\widetilde{C}$ may be defined by
$$
\widetilde{C}(t) = \{x^1(t), \cdots, x^n(t),X^1(t),
\cdots,X^n(t)\}.
$$
where $X^j(t)$ denotes the components of $X$ in a natural basis.
The curve $\widetilde{C}$ is called \emph{horizontal lift} of the
curve $C$ in $M$, if $X$ is  parallel along $C$. When $X$ is the
vector field $\frac{dC}{dt}$ (tangent to $C$), the curve
$\widetilde{C}$ in the tangent bundle is called \emph{the natural
lift} of $C$.

\vskip2mm

Let us consider a curve $C$ in $H_3$ expressed locally by
$x^h=x^h(t)$ and  $Y = y^j(t) \frac{\partial}{\partial x^j}$ a
vector field along $C$. Then, in the tangent bundle $TH_3$, we
define a curve $\widetilde{C}$ by
$$
x^h=x^h(t),\quad y^h=y^h(t),\qquad   h = 1,2,3.
$$

A curve $\gamma(t)=(x^1(t),x^2(t),x^3(t),y^1(t),y^2(t),y^3(t))$ on
$(TH_3,g^s)$ is a geodesic if and only if
\begin{equation*}
\widetilde{\nabla}_{\dot \gamma} \dot \gamma=0
\end{equation*}
where $\widetilde{\nabla}$ is the Levi-Civita connection of the
Sasaki metric on $TH_3,$ and
\begin{equation*}
\dot \gamma=\sum_{i=1}^{3}{\frac{dx^i}{dt}\frac{\delta}{\delta
x^i}+\frac{Dy^i}{dt}\frac{\partial}{\partial y^i}}.
\end{equation*}

Combining the relations above, taking into account the expressions
of Levi-Civita connection for the Sasaki metric, and then
identifying the horizontal and vertical components, we obtain that
$\gamma$ is a geodesic on $(TH_3,g^s)$ if and only if \vskip2mm
\begin{equation*}
\begin{cases}
\frac{D^2x^h}{dt^2}+R^h_{kji}y^k\frac{Dy^j}{dt}\frac{d x^i}{dt}=0\\
\\
\frac{D^2y^h}{dt^2}=0
\end{cases}
\qquad  k,i,j,h = 1,2,3
\end{equation*}
where $R^h_{kji}$ is the curvature of the Heisenberg manifold. We
have denoted
\begin{equation*}\label{Dydt}
\frac{Dy^i}{dt}=\frac{dy^i}{dt}+\Gamma^i_{kj}y^k\frac{dx^j}{dt},
\qquad   i,j,k = 1,2,3.
\end{equation*}

The Lagrangian of the Sasaki metric $g^s$ given by (\ref{Sasaki})
has the expression
\begin{align*}\label{Lagr}
L&=\Big(\frac{dx^1}{dt}\Big)^2+\Big(\frac{dx^2}{dt}\Big)^2+\Big(\frac{dx^3}{dt}+
x^2\frac{dx^1}{dt}-x^1\frac{dx^2}{dt}\Big)^2\\
&+\Big(\frac{Dy^1}{dt}\Big)^2+\Big(\frac{Dy^2}{dt}\Big)^2+
\Big(\frac{Dy^3}{dt}+x^2\frac{Dy^1}{dt}-x^1\frac{Dy^2}{dt}\Big)^2\nonumber
\end{align*}
and the corresponding Euler-Lagrange equations
\[ \frac{d}{dt}\Big(\frac{\delta L}{\delta
\dot{x}^i}\Big)=\frac{\delta L}{\delta x^i},\qquad
\frac{d}{dt}\Big(\frac{\partial L}{\partial
\grave{y}^i}\Big)=\frac{\partial L}{\partial y^i}
\]
where we used the notations $\dot{x}^i=\frac{dx^i}{dt}$ and
$\grave{y}^i=\frac{Dy^i}{dt}$, are

\begin{equation}\label{eq1}
\frac{d}{dt}\Big[\frac{dx^1}{dt}+x^2\Big(\frac{dx^3}{dt}+x^2\frac{dx^1}{dt}-x^1\frac{dx^2}{dt}\Big)\Big]=-\frac{dx^2}{dt}
\Big(\frac{dx^3}{dt}+x^2\frac{dx^1}{dt}-x^1\frac{dx^2}{dt}\Big)
\end{equation}
\begin{equation}\label{eq2}
\frac{d}{dt}\Big[\frac{dx^2}{dt}-x^1(\frac{dx^3}{dt}+x^2\frac{dx^1}{dt}-x^1\frac{dx^2}{dt}\Big)\Big]=\frac{dx^1}{dt}
\Big(\frac{dx^3}{dt}+x^2\frac{dx^1}{dt}-x^1\frac{dx^2}{dt}\Big)
\end{equation}
\begin{equation}\label{eq3}
\frac{d}{dt}\Big(\frac{dx^3}{dt}+x^2\frac{dx^1}{dt}-x^1\frac{dx^2}{dt}\Big)=0
\end{equation}
\begin{equation}
\frac{d}{dt}\Big[\frac{Dy^1}{dt}+x^2\Big(\frac{Dy^3}{dt}+x^2\frac{Dy^1}{dt}-x^1\frac{Dy^2}{dt}\Big)\Big]=0
\end{equation}
\begin{equation}
\frac{d}{dt}\Big[\frac{Dy^2}{dt}-x^1\Big(\frac{Dy^3}{dt}+x^2\frac{Dy^1}{dt}-x^1\frac{Dy^2}{dt}\Big)\Big]=0
\end{equation}
\begin{equation}
\frac{d}{dt}\Big(\frac{Dy^3}{dt}+x^2\frac{Dy^1}{dt}-x^1\frac{Dy^2}{dt}\Big)=0.
\end{equation}
For a geodesic $\gamma:I\rightarrow TH_3,$
$\gamma(t)=(x^1(t),x^2(t),x^3(t),y^1(t),y^2(t),y^3(t))$ which at
the instant zero passes through the origin, with the velocity
$\dot{\gamma}(0)=(u,v,w,l,m,n)$, the Euler-Lagrange equations
become
\begin{equation}\label{Euler-Lagr1}
\frac{dx^3}{dt}+x^2\frac{dx^1}{dt}-x^1\frac{dx^2}{dt}=w
\end{equation}
\begin{equation}\label{Euler-Lagr2}
\frac{d}{dt}\Big(\frac{dx^1}{dt}+x^2w\Big)=-\frac{dx^2}{dt}w
\end{equation}
\begin{equation}\label{Euler-Lagr3}
\frac{d}{dt}\Big(\frac{dx^2}{dt}-x^1w\Big)=\frac{dx^1}{dt}w
\end{equation}
\begin{equation}\label{Euler-Lagr4}
\frac{Dy^3}{dt}+x^2\frac{Dy^1}{dt}-x^1\frac{Dy^2}{dt}=n.
\end{equation}
\begin{equation}\label{Euler-Lagr5}
\frac{Dy^1}{dt}+x^2n=l
\end{equation}
\begin{equation}\label{Euler-Lagr6}
\frac{Dy^2}{dt}-x^1n=m.
\end{equation}

\begin{remark}\label{remark1}
The first three Euler-Lagrange equations, above are satisfied if
and only if the curve $(x^1(t),x^2(t),x^3(t))$ is a geodesic on
the base $(H_3,g)$, which at the moment zero passes through the
point $(0,0,0)$ with the velocity $(u,v,w)$.
\end{remark}

Taking into account that $\dot{\gamma}(0)=(u,v,w,l,m,n)$, from the
equation (\ref{Euler-Lagr3}) we obtain that
$$
\frac{dx^2}{dt}=2x^1w+v
$$
which substituted into (\ref{Euler-Lagr2}), yields the equation
$$
\frac{d^2x^1}{dt^2}+4x^1w^2=0
$$
with the solution
\begin{equation}\label{solx1}
x^1(t)=\frac{v}{2w}\cos(2wt)+\frac{u}{2w}\sin(2wt)-\frac{v}{2w}.
\end{equation}
Analogously, from (\ref{Euler-Lagr2}) and (\ref{Euler-Lagr3}) we
obtain
\begin{equation}\label{solx2}
x^2(t)=-\frac{u}{2w}\cos(2wt)+\frac{v}{2w}\sin(2wt)+\frac{v}{2w}.
\end{equation}

Replacing the solutions (\ref{solx1}) and (\ref{solx2}) into
(\ref{Euler-Lagr1}), we obtain that $x^3$ has the expression
\begin{equation*}
x^3(t)=wt+\frac{u^2+v^2}{2w}t-\frac{u^2+v^2}{2w}\sin(2wt).
\end{equation*}

In the case when $w=0,$ the solutions of the system obtained from
the equations (\ref{Euler-Lagr1}), (\ref{Euler-Lagr2}), and
(\ref{Euler-Lagr3}) have simpler expressions
\begin{equation}\label{rette}
x^1(t)=ut,\ x^2(t)=vt,\ x^3(t)=0.
\end{equation}

We may state now the following result.

\begin{theorem}
The horizontal lift $\widetilde{C}$ and the natural lift
$\widehat{C}$ of a curve $C$ from the Heisenberg manifold $H_3$
are geodesic in the tangent bundle endowed with the Sasaki metric,
$(TH_3,g^s)$   if and only if the curve $C$ is a geodesic in
$(H_3,g),$ and  $\widetilde{C}$, $\widehat{C}$ pass through the
point $(0,0,0,0,0,0),$ such that
$\dot{\widetilde{C}}(0)=\dot{\widehat{C}}(0)=(u,v,w,0,0,0)$.
\end{theorem}
\begin{proof}
If the curve $\widetilde{C}(t)=(C(t),Y(t))$ is the horizontal lift
to $TH_3$ of the curve $C(t)=(x^1(t),x^2(t),x^3(t))$ from $H_3$,
then $Y$ is a parallel vector field along $C$, i.e.
$$
\frac{Dy^1}{dt}=\frac{Dy^2}{dt}=\frac{Dy^3}{dt}=0
$$
and in this case, the last three Euler-Lagrange equations,
(\ref{Euler-Lagr4}), (\ref{Euler-Lagr5}), and (\ref{Euler-Lagr6})
reduce to $l=m=n=0.$

If the curve $\widehat{C}(t) = (C(t),Y(t))$ on $TH_3$ is the
natural lift of the curve $C(t)=(x^1(t),x^2(t),x^3(t))$ from
$H_3$, then $Y$ is the tangent vector field to $C$, i.e.
$$
y^h=\frac{dx^h}{dt},\qquad  h = 1,2,3
$$
from which we obtain that the covariant derivative of $Y$ has
the expression
\begin{equation}\label{Dy}
\frac{Dy^h}{dt}=\frac{d^2x^h}{dt^2}+
\Gamma^h_{ij}\frac{dx^i}{dt}\frac{dx^j}{dt},\quad \forall
i,j,h=\overline{1,3}.
\end{equation}

Taking into account Remark \ref{remark1}, it follows that the
curve $C(t)=(x^1(t),$ $x^2(t),x^3(t))$ is a geodesic on the
Heisenberg manifold $(H_3,g)$, and then $\frac{Dy^h}{dt}$
expressed by (\ref{Dy}) vanishes, and the equations
(\ref{Euler-Lagr4}), (\ref{Euler-Lagr5}), (\ref{Euler-Lagr6})
reduce again to $l=m=n=0.$ Thus the theorem is proved.
\end{proof}

In a more general context, when we search some examples of
geodesics on $(TH_3,g^s)$, which are not horizontal or natural
lifts of the curves from the base manifold $(H_3,g)$, we obtain,
by taking into account the expressions of the Christoffel symbols
constructed with the Heisenberg metric, that the last three
Euler-Lagrange equations, (\ref{Euler-Lagr4}),
(\ref{Euler-Lagr5}), (\ref{Euler-Lagr6}) get the forms
\begin{align*}
\frac{dy^1}{dt} &+
x^2y^2\frac{dx^1}{dt}+(x^2y^1-2x^1y^2+y^3)\frac{dx^2}{dt}+y^2\frac{dx^3}{dt}=l-vnt
\\
\frac{dy^2}{dt}&+(-2x^2y^1+x^1y^2-y^3)\frac{dx^1}{dt}+x^1y^1\frac{dx^2}{dt}-y^1\frac{dx^3}{dt}=m+unt
\\
\frac{dy^3}{dt}&+[-2x^1x^2y^1+(1+(x^1)^2-(x^2)^2)y^2-x^1y^3]\frac{dx^1}{dt}
\\
&+ [(1+(x^1)^2-(x^2)^2)y^1+2x^1x^2y^2-x^2y^3]\frac{dx^2}{dt}
\nonumber \\
&- (x^1y^1+x^2y^2)\frac{dx^3}{dt}+x^2(l-x^2n)-x^1(m+x^1n)=n. \nonumber
\end{align*}

In the case when the curve on the base manifold is a geodesic
given by (\ref{rette}), then the above equations become
\begin{align*}
\frac{dy^1}{dt}&+ tvy^2u+(tvy^1-2tuy^2+y^3)v=l-vnt
\\
\frac{dy^2}{dt}&+(-2tvy^1+tuy^2-y^3)u+tuy^1v=m+unt
\\
\frac{dy^3}{dt}&+[-2t^2uvy^1+(1+t^2u^2-t^2v^2)y^2-tuy^3]u
\\
&+ [(1+t^2u^2-t^2v^2)y^1+2t^2uvy^2-tvy^3]v
\nonumber \\
&+tv(l-vnt)-tu(m+unt)=n \nonumber
\end{align*}
i.e. we have the system

\begin{equation}\label{sist}
\begin{cases}
\frac{dy_1}{dt}+ tv^2y_1-tuvy_2+vy_3=l-vnt\\
\frac{dy_2}{dt}-tuvy_1+tu^2y_2-uy_3=m+unt\\
\frac{dy_3}{dt}+ v(1-
t^2u^2-t^2v^2)y_1+u(1+t^2u^2+t^2v^2)y_2-t(u^2+v^2)y_3 +\\
 \qquad +tv(l-vnt)-tu(m+unt)=n.
\end{cases}
\end{equation}

\begin{remark} \rm If $u=v=0,$ we obtain the following
particular solution of the system (\ref{sist}):
$$
y_1=lt,\ y_2=mt,\ y_3=nt.
$$
\end{remark}

Taking this remark into account, we may prove:

\begin{theorem}\label{th}
If the curve from the Heisenberg manifold reduces to the origin
point $(0,0,0),$ then the geodesics from the tangent bundle with
the Sasaki metric $(TH_3,g^s)$ are the curves $\gamma$ passing
through the origin $(0,0,0,0,0,0)$ with velocity
$\dot{\gamma}(0)=(0,0,0,l,m,n)$, namely
$\gamma(t)=(0,0,0,lt,mt,nt).$
\end{theorem}

Yanno and Ishihara proved that if a geodesic lies in a fiber of
the tangent bundle $(TM,g^s)$ of an $n$-dimensional Riemannian
manifold $(M,g)$, given by $x^h=c^h, \forall h=\overline{1,n}$,
where $c^h$ is a real constant, then the geodesic is expressed by
linear equations $x^h=c^h,\ y^h=a^ht+b^h,$ with respect to the
induced coordinates $\{x^h,y^h\}_{h=\overline{1,n}},$ where
$a^h,b^h,c^h$ are constants. In the case of the tangent bundle
$(TH_3,g^s)$ this result reduces to Theorem \ref{th}, since when
$x^i$ are constants, the expressions (\ref{Dydt}) of
$\frac{Dy^i}{dt}$ become $\frac{Dy^i}{dt}=\frac{dy^i}{dt},
i=\overline{1,3}$, and the Euler-Lagrange equations
(\ref{Euler-Lagr4}) - (\ref{Euler-Lagr6}) with the initial
conditions
$\frac{dy^1}{dt}(0)=l,\frac{dy^2}{dt}(0)=m,\frac{dy^3}{dt}(0)=n$
lead to $x^i=0,i=\overline{1,3}$.

\begin{remark}\label{remark} If $m=n=0$, a particular solution of the
system (\ref{sist}) is of the form
\begin{equation*}
y_1=lt,\ y_2=0,\ y_3=-lvt^2.
\end{equation*}
\end{remark}

Taking into account Remark \ref{remark} and the solution
(\ref{rette}) of the system obtained from the equations
(\ref{eq1}), (\ref{eq2}), (\ref{eq3}), we may formulate:

\begin{theorem}
One of the geodesics from the tangent bundle with the Sasaki
metric $(TH_3,g^s)$, is a curve $\widetilde{\gamma}:I\rightarrow
TH_3,$ which at the moment zero passes through the point
$\widetilde{\gamma}(0)=(0,0,0,0,0,0),$ with the property
$\dot{\widetilde{\gamma}}(0)=(u,v,0,l,0,0)$, namely the curve
$\widetilde{\gamma}(t)=(ut,vt,0,lt,0,-lvt^2)$.
\end{theorem}

\bf{Acknowledgements.} \rm The authors would like to express their
gratitude to professors R. Caddeo, V. Oproiu and M. I. Munteanu,
for the useful hints and advices throughout this work. 
\strut\hfill Al.I. Cuza University of Ia\c si,\\
\strut\hfill Faculty of Mathematics\\
\strut\hfill Bd. Carol I Nr. 11, 700 506  Ia\c si, ROM\^{A}NIA\\
\vspace{-2mm}
\strut\hfill simonadruta@yahoo.com \\

\strut\hfill Universit\`a degli Studi di Cagliari,\\
\strut\hfill Dipartimento di Matematica e Informatica\,\\
\strut\hfill Via Ospedale 72, 09124 Cagliari, ITALIA\\
\vspace{-2mm}
\strut\hfill piu@unica.it
\bigskip

\end{document}